\documentclass[12pt, reqno]{amsart}
\usepackage{amsmath, amsthm, amscd, amsfonts, amssymb, graphicx, color}
\usepackage[bookmarksnumbered, colorlinks, plainpages]{hyperref}

\newtheorem{theorem}{Theorem}[section]
\newtheorem{lemma}[theorem]{Lemma}
\newtheorem{proposition}[theorem]{Proposition}
\newtheorem{corollary}[theorem]{Corollary}
\theoremstyle{definition}
\newtheorem{definition}[theorem]{Definition}
\newtheorem{example}[theorem]{Example}

\theoremstyle{remark}
\newtheorem{remark}[theorem]{Remark}
\numberwithin{equation}{section}

\begin{document}

\title[$\mathcal A$-2-Frames  in $\mathcal A$-2-inner product spaces]
{$\mathcal A$-2-Frames  in $\mathcal A$-2-inner product spaces}

\author[ B. Mohebbi Najmabadi and  T.L. Shateri  ]{B. Mohebbi Najmabadi and  T.L. Shateri }
\address{Behrooz Mohebbi Najmabadi \\ Phd student of Mathematics, Hakim Sabzevari University, Sabzevar, IRAN }
\email{\rm behrozmohebbi1351@gmail.com; behruz\_mohebbi@yahoo.com}
\address{Tayebe Lal Shateri \\ Department of Mathematics and Computer
Sciences, Hakim Sabzevari University, Sabzevar, P.O. Box 397, IRAN}
\email{ \rm t.shateri@hsu.ac.ir; t.shateri@gmail.com}

\thanks{*The corresponding author:
t.shateri@hsu.ac.ir; t.shateri@gmail.com  (Tayebe Lal Shateri)}
\date{September, 14, 2019}

\maketitle

\begin{abstract}
\normalsize 
Certain results about  frames are extended for the new frames in Hilbert $C^*$-modules. In this paper, we introduce the notion of  $\mathcal A$-2-frames in $\mathcal A$-2-inner product spaces and give some characterizations for these frames. Then we define the tensor product of $\mathcal A$-2-frames and prove some results for it.
\end{abstract}
\section{Introduction and preliminaries}
 Frames were first introduced in 1952 by Duffin and Schaeffer \cite{JA}. In 2000, Frank-Larson \cite{MD} introduced the notion of frames in Hilbert $ C^*$-modules as a generalization of frames in Hilbert spaces  and Jing \cite{WJ} continued  to consider them. It is well known that Hilbert $C^*$-modules are  generalizations of Hilbert spaces by allowing the inner product to take values in a $C^*$-algebra rather than in the field of complex numbers. The theory of $2$-inner product spaces as well as an extensive list of related references can be found in \cite{SG,YPL}. The concept of 2-frames for 2-inner product spaces was introduced by A. Arefijamal and Ghadir Sadeghi \cite{AS} and described some fundamental properties of them. 
Recently T. Mehdiabad and A. Nazari \cite{TMN} introduced the $\mathcal A$-2-inner product space and investigate some inequalities in these spaces. The authors \cite{MS} defined a 2-inner product that takes values in a locally $C^*$-algebra and studied some properties of it.\\   
In this paper, we introduce an $\mathcal A$-2-frame in the $\mathcal A$-2-inner product space and describe some fundamental properties of them. The tensor product of $\mathcal A$-2-frames in the $\mathcal A$-2-inner product space is introduced. It is shown that the tensor product of two $\mathcal A$-2-frames is an $\mathcal A$-2-frame for the tensor product of $\mathcal A$-2-inner product space. Also, we investigate tensor products of $\mathcal A$-2-frames.
From now,  $\mathcal A$ denotes a $C^*$-algebra.
\begin{definition} 
 A pre-Hilbert $\mathcal A$-module is a complex vector space $E$ which is also a left $\mathcal A$-module, compatible with the Complex algebra structure,  equipped with an $\mathcal A$-valued inner product \\$\langle .,.\rangle:E\times E\to \mathcal A $  which is $\mathbb C$-linear and $\mathcal A$-linear in its second variable and satisfies the following relations\\
$(I_1)$  $\langle x ,x \rangle\geq0$  for  every $ x\in E$,\\
$(I_2)$  $\langle x ,y \rangle =\langle y ,x \rangle^*$ for  every $ x ,y \in E$,\\
$(I_3)$  $\langle x ,x \rangle =0 $ if  and  only  if $  x=0 $,\\
$(I_4)$  $\langle ax ,by \rangle =a^*\langle x ,y \rangle b $ for  every $  x ,y \in E$  and $ a ,b \in \mathcal A$,\\
$(I_5)$  $\langle x ,\alpha y +\beta z\rangle =\alpha\langle x ,y\rangle +\beta\langle x ,z\rangle $  for  every $x ,y ,z \in E$ and $  \alpha ,\beta \in \mathbb{C}$.
\end{definition}
\begin{example}
Let $l^{2}(\mathcal A)$ be the set of all sequences $\{a_{n}\}_{n\in\mathbb{N}}$ of elements of a  $C^*$-algebra $\mathcal A$ such that the series $\sum_{n\in \mathbb{N}}a_{n}a^*_{n}$ is convergent in $\mathcal A$. Then $l^{2}(\mathcal A)$ is a Hilbert $\mathcal A$-module with respect to the pointwise operations and inner product defined
\begin{equation*}
\langle \{a_n\}_{n\in \mathbb N} , \{b_n\}_{n\in \mathbb N}\rangle =\sum_{n\in  \mathbb N}a_n b^*_n.
\end{equation*}
\end{example}
\begin{definition}
Let $E$ be a left $\mathcal A$-module, an $\mathcal A$-combination of $x_{1},x_{2},...,x_{n}$ in $E$ is written as follows
\begin{equation*}
 \sum_{i=1}^{n} a_{i}x_{i}=a_{1}x_{1}+a_{2}x_{2}+...+a_{n}x_{n}  \quad {(a_{i} \in \mathcal A)}.
\end{equation*}
 $x_{1},x_{2},...,x_{n}$ are called $\mathcal A$-independent if the equation $a_{1}x_{1}+a_{2}x_{2}+...+a_{n}x_{n}=0$ has exactly one solution, namely $a_{1}=a_{2}=...=a_{n}=0$, otherwise, we say that  $x_{1},x_{2},...,x_{n}$ are $\mathcal A$-dependent.\\The maximum number of elements  in $E$ that are $\mathcal A$-independent is called the $\mathcal A$-rank of $E$.
\end{definition}
\begin{definition}
Let $\mathcal A$ be a C*-algebra and $ E$ be a linear space by$\mathcal A$-rank greater than 1,  which is also a left $\mathcal A$-module.  We define a function $ \langle . , .|  .\rangle  :E\times E\times E \to \mathcal A$ satisfies  the  following  properties\\
$(T_1)$ $ \langle x,x|y\rangle  =0$, if and only if $x=ay$  for $ a\in\mathcal A $  \\
$(T_2)$ $ \langle x , x | y\rangle  \geq 0 $ for all $x ,y \in E$ \\
$(T_3)$ $ \langle x ,x | y\rangle  =\langle y , y | x\rangle $ for all $x ,y \in E $ \\
$(T_4)$ $ \langle x ,y | z\rangle  =\langle y , x | z\rangle ^*$ for all $x ,y ,z \in E $\\
$(T_5)$ $ \langle ax ,by | z\rangle  = a\langle x , y | z\rangle  b^*$ for all $x ,y ,z \in E$  and $ a, b \in \mathcal A$ \\
$(T_6)$ $ \langle x ,\alpha y | z \rangle = \overline{\alpha}\langle x , y | z\rangle $ for all$x ,y \in E $ and $\alpha\in \mathbb{C} $\\
$(T_7)$ $ \langle x+y ,z | w\rangle  =\langle x , z | w\rangle  +\langle y , z | w\rangle $ for all $x ,y , z ,w \in E  .$\\
Then the function $ \langle . , .|  .\rangle $ is called $ \mathcal A$-2- inner product  and $ ( E , \langle . , .|  .\rangle  ) $ is called  $\mathcal A$-2-inner product space. 
\end{definition}
\begin{example}\cite{MS}
Let $\mathcal  A $ be a commutative $C^*$-algebra and $E$ be a  pre-Hilbert $\mathcal A$-module  with inner product $\langle . , . \rangle,$  define
$ \langle . , .|  .\rangle $ :$ E\times E\times E \to \mathcal A$ by \\
\begin{align*}
 (x ,y ,z)\longmapsto \langle x , y | z \rangle  =\langle x , y\rangle \langle z ,z\rangle  -\langle x ,z\rangle\langle z ,y\rangle
\end{align*}
Then  $ ( E , \langle . , . |  .\rangle ) $ is a  $\mathcal A$-2- inner product space.
\end{example}
\begin{theorem}
Let $(E,\langle . , . | .\rangle )$ be an $ \mathcal A$-2-inner product space on a  commutative $C^*$-algebra $\mathcal A$. Then the following inequality holds,\\
\begin{align*}
|\langle x , y | z\rangle |^{2}=\langle x , y | z\rangle \langle x , y | z\rangle ^{*}\leq \langle x , x | z\rangle \langle y , y | z\rangle \quad {(x,y,z\in E).}
\end{align*}
\begin{proof}
For $\lambda \in \mathcal A$ we have 
\begin{align*}
0\leq\langle \lambda x-y ,\lambda x-y|z\rangle &=\langle \lambda x,\lambda x|z\rangle -\langle\lambda x,y| z\rangle -\langle y ,\lambda x|z\rangle +\langle y,y|z\rangle \\&=\lambda^*\langle x,x|z\rangle \lambda-\lambda^*\langle x,y| z\rangle -\langle y,x|z\rangle \lambda+\langle y,y|z\rangle.
\end{align*}
Take $\lambda=\langle x,y|z\rangle (\langle x,x|z\rangle+\varepsilon e)^{-1}$ then\\
\begin{align*}
0&\leq \langle y,x|z\rangle (\langle x,x|z\rangle +\varepsilon e)^{-1}\langle x,x|z\rangle \langle x,y|z\rangle(\langle x,x|z\rangle +\varepsilon e)^{-1}\\&-\langle y,x|z\rangle (\langle x,x|z\rangle +\varepsilon e)^{-1}\langle x,y|z\rangle -\langle y,x|z\rangle \langle x,y|z\rangle (\langle x,x|z\rangle +\varepsilon e)^{-1}+\langle y,y|z\rangle,
\end{align*} 
hence,
$2\langle y,x|z\rangle \langle x,y|z\rangle \leq (\langle x,x|z\rangle +\varepsilon e)^{-1}\langle x,x|z\rangle \langle y,x|z\rangle \langle x,y|z\rangle \\+\langle y,y|z\rangle (\langle x,x|z\rangle +\varepsilon e)\leq (\langle x,x|z\rangle+\varepsilon e)^{-1}(\langle x,x|z\rangle +\varepsilon e)\langle y,x|z\rangle \langle x,y|z\rangle +\langle y,y|z\rangle (\langle x,x|z\rangle +\varepsilon e)$
then by $\varepsilon\rightarrow 0$ inequality holds.
\end{proof}
\end{theorem}
\begin{definition}\cite{TMN}
Let $E$ be a real vector space that $\mathcal A$-rank is greater than 1 and\\ $p:E\times E\to\mathbb{ R}$ be a function such that\\
$(1)$ $p(x , y) =0$  if and only if $x, y \in E $ are linearly $\mathcal A$ - dependent, \\
$(2)$ $p(x , y)=p(y , x)$ for every $x , y \in E,$\\
$(3)$ $p(\alpha x, y)=|\alpha|p(x ,y )$, for every  $x , y\in E$ and for every $\alpha\in\mathbb{C},$\\
$(4)$ $p(x+y , z )\leq p(x , z)+p(y , z)$, for every $x , y ,z \in E.$\\
$(5)$ $P(ax , y )\leq||a||p(x , y )$ ,  for every $x, y \in E$ and $a\in \mathcal A$,
The function $p$ is called an $\mathcal A$-2-norm.
\end{definition}
It follows from theorem 1.6 that
\begin{corollary}\cite{TMN}
Let $E$ be an  $\mathcal A$-2- inner product space, For $x,z\in E$ we define $ p(x , z ) = \sqrt{\Arrowvert\Big(\langle x , x | z\rangle \Big)\Arrowvert}$. Then $ ||\langle x , y | z\rangle||\leq p(x, z )p(y, z )$.
\end{corollary}

In the following theorem, we investigate some properties of an $ \mathcal A$-2-norm.
\begin{theorem}
Let $(E,\langle . , . | .\rangle )$ be an $\mathcal A$-2-inner product space and $p$ be an  $\mathcal A$-2-norm,  then \\
$(1)$ $p(x, y) = \sup\left\{\Arrowvert\langle x, z|y \rangle\Arrowvert;p(z, y)= 1\right\}$.\\
$(2)$ $p(x,y+ax)=p(x,y)$ for $a\in \mathcal A$.
\begin{proof}
$(1)$ By the Cauchy-schwarz inequality we observe that
$\Arrowvert\langle x, z|y \rangle\Arrowvert\leq p(x, y)p(z, y)\leq p(x,y),$
for every $z\in E$ such that $p(z, y)\leq 1$. 
Moreover if $z=\dfrac{x}{p(x,y)}$ then $p(z,y)=1$ and therefore  $\Arrowvert\langle x, z| y\rangle\Arrowvert=p(x,y)$.
\end{proof}
\end{theorem}
Let $E$ be an $\mathcal  A$-2- inner product space. A sequence $\{a_{n}\}_{n\in \mathbb{N}}$ of $E$ is said to be convergent if there exists an element $a\in E$ such that $\lim_{n\rightarrow\infty}p(a_{n}-a,x)=0$,  for all $x\in E$. Similarly, we can define a Cauchy sequence in $E$. An $\mathcal A$-2- inner product space $E$ is called  an  $\mathcal A$-2- Hilbert  space if it is complete.\\
Now, we give the notion of a frame on a Hilbert $\mathcal A$-module which is defined in \cite[definition 3.1]{WJ}.
\begin{definition}\cite{WJ}
Let $\mathcal{A}$ be an unital  $C^*$-algebra and $ E$ be a  Hilbert $ \mathcal{A}$-module.  The sequence $\{x_{j}\in~ E | j\in J\subseteq\mathbb{N}\}$ is called  a frame for $E$ if there exist two  positive elements $A$ and $B$ in real numbers such that
\begin{align*}
A\langle x,x\rangle \leq\sum_{j\in J}\langle x,x_{j}\rangle \langle x_{j},x\rangle \leq B\langle x,x\rangle  \qquad( x\in  E).
\end{align*}
The frame $\{x_{j}\}$ is said to be tight frame if $A=B$, and said to be Parseval if $A=B=1$.\\
The operator $T:E\rightarrow {l}^2(\mathcal{A})$ defined by
\begin{align*}
Tx=\{\langle x,x_{j}\rangle \}_{j\in J}
\end{align*}
is called the analysis operator. The adjoint operator $T^*: {l}^2(\mathcal{A}) \rightarrow E$ is given by
\begin{align*}
T^*\{c_{j}\}_{j\in J}=\sum_{j\in J}c_{j}x_{j}
\end{align*}
 is called the pre-frame operator or the synthesis operator. By composing $T$ and $T^*$, we obtain the frame operator $ S:E\rightarrow E,$ by
\begin{align*}
S=T^*Tx=\sum_{j\in J}\langle x,x_{j}\rangle x_{j}.
\end{align*}
Also from this equation, we have
\begin{align*}
x=\sum_{j\in J}\langle x,S^{-1}x_{j}\rangle x_{j}  \qquad( x\in  E).
\end{align*}
\end{definition}
Now  we are ready to define an $\mathcal A$-2-frame on an  $\mathcal A$-2- Hilbert space.
\section{ $\mathcal A$-2-frames}
In this section we define $\mathcal A$-2-frames  on $\mathcal A$-2- Hilbert spaces, and  we give some results about them.
\begin{definition}
Let $(E,\langle ., . | .\rangle )$ be an $\mathcal A$-2- Hilbert space and $\xi \in E$. A sequence  $\{a_{i}\}_{i\in \mathbb{N}}$ of $E$ is called an $\mathcal A$-2-frame (associated to $\xi$) if there exist  positive real  numbers $A$ and $B$ such that
\begin{align}\label{2.1}
A\langle x,x|\xi\rangle \leq\sum_{i\in \mathbb{N}}\langle x , x_{i} |\xi\rangle \langle x_{i}, x |\xi\rangle \leq B\langle x ,x |\xi \rangle \qquad(x\in  E ).
\end{align}
A sequence satisfying the upper $\mathcal A$-2-frame condition is called an $\mathcal A$-2-Bessel sequence, and every $x_{i}$ is $A$- independent to $\xi$.
\end{definition}
\begin{proposition}
Let $\mathcal A$ be a commutative and  $ (E,\langle . , .\rangle )$ be a Hilbert $\mathcal A$-module and $\{x_{i}\}_{i\in \mathbb{N}}$ be a frame for $E$. Then  for invertible element $\langle \xi, \xi \rangle $,  it is  an $\mathcal A$-2-frame with the standard $\mathcal A$-2-inner product.
\begin{proof}
\begin{align*}
\sum_{i\in \mathbb{N}}\langle x , x_{i} |\xi\rangle \langle x_{i}, x |\xi\rangle &=\sum_{i\in \mathbb{N}}\langle x\langle \xi, \xi \rangle -\langle \xi, x\rangle \xi , x_{i}\rangle \langle x_{i} ,x\langle \xi, \xi \rangle -\langle \xi, x\rangle \xi\rangle \\&\leq B\langle x\langle \xi, \xi \rangle -\langle \xi, x\rangle \xi,x\langle \xi, \xi \rangle -\langle \xi, x\rangle \xi\rangle \\&\leq B\langle \xi, \xi \rangle \Big(\langle x , x\rangle \langle \xi , \xi\rangle -\langle x , \xi\rangle \langle \xi , x\rangle \Big)\\&=B\langle \xi, \xi \rangle \Big(\langle x , x |\xi\rangle \Big)\leq ||B\langle \xi, \xi \rangle ||\Big(\langle x , x |\xi\rangle \Big)
\end{align*}
Take $D=||B\langle \xi , \xi\rangle ||$, the argument for lower bound is similar.
\end{proof}
\end{proposition}

In the following proposition, $E$ is a  Hilbert $\mathcal A$-module in which every closed submodule is orthogonally complemented and $\langle \xi, \xi \rangle $ is invertible and $L_{\xi}$ is the subspace generated with $\xi$.
\begin{proposition}
Let $\mathcal A$ be a commutative and $ (E,\langle . , .\rangle )$ be a Hilbert $\mathcal A$-module and $\xi\in E$. Every $\mathcal A$-2-frame associated with $\xi$ is a frame for $L^{\bot}_{\xi}$.
\begin{proof}
\begin{align*}
\\&\sum_{i\in \mathbb{N}}\langle x\langle \xi, \xi \rangle -\langle \xi, x\rangle \xi , x_{i}\rangle \langle x_{i} ,x\langle \xi \rangle -\langle \xi, x\rangle \xi\rangle \\&\leq B\langle x\langle \xi, \xi \rangle -\langle \xi, x\rangle \xi,x\langle \xi, \xi \rangle -\langle \xi, x\rangle \xi\rangle \\&\leq B\langle \xi, \xi \rangle \Big(\langle x , x\rangle \langle \xi , \xi\rangle -\langle x , \xi\rangle \langle \xi , x\rangle \Big)
\end{align*}Then
\begin{align*}
A\langle x , x \rangle \langle \xi, \xi\rangle ^{2}\leq\langle \xi , \xi\rangle ^{2}\sum_{i\in \mathbb{N}}\langle x , x_{i}\rangle \langle x_{i} , x\rangle \leq B\langle \xi , \xi \rangle ^{2}\langle x , x\rangle   \qquad(x\in L^{\bot}_{\xi})
\end{align*}
Since $\langle \xi , \xi \rangle $ is invertible, the  proof is completed.
\end{proof}
\end{proposition}
Let $(E,\langle ., . | .\rangle )$ be an $\mathcal A$-2-Hilbert space and $L_{\xi}$ be the subspace generated with $\xi$ for a fix element $\xi$ in $E$. Denote by $\mathcal M_{\xi}$ the algebraic complement of $L_{\xi}$ in $E$. So  $L_{\xi}\oplus\mathcal M_{\xi}=E$.
We define  the semi-inner product $ \langle . , .\rangle _{\xi}$ on $E$ as following
\begin{align*}
 \langle x, z\rangle _{\xi}=\langle x , z |\xi\rangle .
\end{align*}
This semi-inner product induces an inner product on the quotient space $E/L_{\xi}$ as
\begin{align*}
\langle x+L_{\xi}, z+L_{\xi}\rangle _{\xi}=\langle x , z \rangle_{\xi}.  \qquad(z,  x\in  E).
\end{align*}
By identifying  $E/L_{\xi}$ with $\mathcal M_{\xi}$ in an obvious way, we obtain an inner product on $\mathcal M_{\xi}$.\\
Now if $\{x_{i}\}_{i\in \mathbb{N}}\subseteq E$  is an $\mathcal A$-2- frame associated with $\xi$ with bounds $A$ and $B$, we can rewrite(2.1) as 
\begin{align*}
A\langle x , x \rangle _{\xi}\leq\sum_{i\in \mathbb{ N}}\langle x , x_{i}\rangle \langle x_{i},x\rangle \leq B\langle x , x\rangle _{\xi}     \qquad(x\in\mathcal {M_{\xi}}).
\end{align*}
That is, $\{x_{i}\}_{i\in \mathbb{N}}$ is a frame for  $\mathcal M_{\xi}$. Let $ E_{\xi}$ be the completion of the inner product space  $ E_{\xi}$, then the sequence  $\{x_{i}\}_{i\in \mathbb{N}}$ is also a frame for $ E_{\xi}$. To summarize, we have the following theorem.
\begin{theorem}
Let $(E,\langle ., . | .\rangle )$ be an $\mathcal A$-2-Hilbert space. Then  $\{x_{i}\}_{i\in \mathbb{N}}\subseteq E$ is  an $\mathcal A$-2-frame associated with $\xi$ if and only if it is a frame for the Hilbert space  $ E_{\xi}$.
\end{theorem}
\begin{lemma}
Let  $\{x_{i}\}_{i\in \mathbb{N}}$ be an $\mathcal A$-2-Bessel sequence in $E$. Then the $\mathcal A$-2-pre frame operator $T:l^2(\mathcal A)\rightarrow E_{\xi}$
defined by
\begin{align*}
T_{\xi}\{c_{i}\}=\sum_{i\in \mathbb{N}}c_{i}x_{i}
\end{align*}
is well-defined and bounded.
\begin{proof}
\begin{align*}
\| \sum_{i=1}^{n} c_{i}x_{i}- \sum_{i=1}^{m} c_{i}x_{i} , \xi \| ^{2}&=sup\{\| \langle \sum_{i=m+1}^{n} c_{i}x_{i} , y | \xi \rangle  \|^{2}  , y\in E_{\xi} ,\|  y , \xi  \|=1\}\\&\leq sup\{\| \langle  \sum_{i=m+1}^{n} \langle x_{i} , y | \xi \rangle 
\langle y , x_{i} | \xi \rangle  \|  , y\in E_{\xi} ,\|  y , \xi  \|=1\} \| \sum_{i=m+1}^{n} c_{i}c^*_{i} \|\\& \leq \| \sum_{i=m+1}^{n} c_{i}c^*_{i} \| \| B \| .
\end{align*}
Since $B$ is the upper bound of $\{x_{i}\}$ this implies that $\sum_{i=1}^{\infty}c_{i}x_{i}$ is well defined as an element of $E_{\xi}$. Moreover if $\{c_{i}\}$ is a sequence in $ l^2(\mathcal A)$, then an argument as above shows that
\begin{align*}
\| T_{\xi}(\{c_{i}\}) \| \leq \sqrt{\| B \| } \| \sum_{i \in \mathbb{N}}c_{i}c^*_{i} \|
\end{align*}
In particular, $\| T_{\xi}\| \leq  \sqrt{\| B \| }$.
\end{proof}
\end{lemma}
We have $\langle x , T(c_{j}) | \xi \rangle=\langle x , T(c_{j}) \rangle_{\xi}=\langle x ,\sum c_{j}x_{j} \rangle_{\xi}=\sum \langle x ,x_{j} \rangle_{\xi}c_{j}^*=\langle {\langle x, x_{j} \rangle_{\xi}} , {c_{j}} \rangle$
Next, we can compute $T^*_{\xi}$, the adjoint of  $T_{\xi}$ as
\begin{align*}
T^*_{\xi}: E_{\xi}\rightarrow l^2(\mathcal A);    \qquad T^*_{\xi}x=\{\langle x , x_{i}|\xi\rangle \}_{i\in \mathbb{N}}.
\end{align*}
 $T^*_{\xi}$ is well-defined and bounded, because
\begin{align*}
\| T^*_{\xi}(x) \|^{2}=\| \{\langle x , x_{i} ,|\xi \rangle \}_{i \in \mathbb{N}} \|^{2}=\| \sum_{i \in \mathbb{N}}\langle x , x_{i} ,|\xi \rangle \langle x_{i} , x | \xi \rangle  \| \leq \|B\| \|x , \xi \|
\end{align*}
That implies $\| T^*_{\xi} \| \leq \sqrt {\| B \|}$.
\begin{definition}\label{2.2}
Let $ \{x_{i}\}_{i\in \mathbb{N}}$ be an $\mathcal A$-2-frame  associated to $\xi$   with bounds $A$ and $B$ in an $\mathcal A$-2- Hilbert space $E$. The operator $S_{\xi}:E_{\xi}\rightarrow E_{\xi}$ defined by
\begin{align}
S_{\xi}x=\sum_{i\in \mathbb{N}}\langle x , x_{i} |\xi\rangle x_{i}.
\end{align}
 is called the $\mathcal A$-2-frame operator for $ \{x_{i}\}_{i\in \mathbb{N}}$.
\end{definition}
In the next theorem, we investigate some properties of $S_{\xi}.$
\begin{theorem}
Let $ \{x_{i}\}_{i\in \mathbb{N}}$ be an $\mathcal A$-2-frame  associated to $\xi$ for an $\mathcal A$-2-Hilbert space $(E ,\langle . ,. | \rangle )$ with $\mathcal A$-2-frame operator $ S_{\xi}$. Then $ S_{\xi}$ is bounded,  invertible, self-adjoint, and positive.
\begin{proof}
It is clear that $ S_{\xi}=T_{\xi}T^*_{\xi}$ is self adjoint and
\begin{align*}
\| S_{\xi} \|=\| T_{\xi}T^*_{\xi} \|=\| T_{\xi} \|^{2} \leq \| B \|.
\end{align*}
We can conclude the boundedness of $S_{\xi}$ directly
\begin{align*}
\| S_{\xi}(x) , \xi \|^{2}&=sup\{\| \langle  S_{\xi}(x) , y | \xi \rangle  \|^{2}  , y\in E_{\xi} ,\|  y , \xi  \|=1\}\\&= sup\{\| \langle \sum_{i=1}^{\infty}\langle x , x_{i} | \xi \rangle x_{i} , y | \xi \rangle  \|^{2} , y\in E_{\xi} ,\|  y , \xi  \|=1\}\\& \leq sup\{\| \sum_{i=1}^{\infty}\langle x , x_{i} | \xi \rangle \langle x_{i} , x | \xi \rangle \| \| \sum_{i=1}^{\infty}\langle y , x_{i} | \xi \rangle \langle x_{i} , y | \xi \rangle \|  ,\|  y , \xi  \|=1 \}\\& \leq \| B \|^{2} \| \langle  x , x | \xi \rangle \|
\end{align*}
The inequality (\ref{2.1}) means that 
\begin{align*}
A\langle x,x \rangle _{\xi}\leq \langle S_{\xi}(x), x \rangle _{\xi}\leq B\langle x ,x \rangle _{\xi} 
\end{align*}
$ S_{\xi} $ is a positive element in the set of all bounded operators on the Hilbert space $E_{\xi}$.
\end{proof}
\end{theorem}
By the definition of $S_{\xi}$ we get the following results.
\begin{corollary}
Let $ \{x_{i}\}_{i\in \mathbb{N}}$ be an $\mathcal A$-2-frame in an $\mathcal A$-2-Hilbert space $(E ,\langle . ,. | \rangle )$ with frame operator  $ S_{\xi}$. Then each $x\in E_{\xi}$ has an expansion of  the following
\begin{align*}
x= S S_{\xi}^{-1}x=\sum_{i\in \mathbb{N}}\langle  S_{\xi}^{-1}x , x_{i} |\xi\rangle x_{i}.
\end{align*}
\end{corollary}
\begin{corollary}
Let $\xi$ and $\eta$ be $\mathcal A$-independent and $ \{x_{i}\}_{i\in \mathbb{N}}$ be an $\mathcal A$-2-frame associated with $\xi$ and $\eta$,  and for $x\in E_{\xi}\cap E_{\eta}$, the operators $ S_{\xi}, E_{\eta}:E_{\xi}\cap E_{\eta}\rightarrow E_{\xi}\cap E_{\eta}$ defined by,
\begin{align*}
S_{\xi}x=\sum_{i\in \mathbb{N}}\langle x , x_{i} |\xi\rangle x_{i},\\
S_{\eta}x=\sum_{i\in \mathbb{N}}\langle x , x_{i} |\eta\rangle x_{i}.
\end{align*}
Then we have $\langle S_{\eta}x, x | \xi \rangle =\langle S_{\xi} x, x | \eta \rangle^*$.
\begin{proof}
\begin{align*}
\langle S_{\eta}x, x | \xi \rangle &=\langle \sum_{i\in \mathbb{N}}\langle x , x_{i} |\eta\rangle x_{i}, x | \xi \rangle\\& =\sum_{i\in \mathbb{N}}\langle x , x_{i} |\eta\rangle \langle x_{i}, x | \xi \rangle\\&=(\sum_{i\in \mathbb{N}}\langle x , x_{i} |\xi \rangle \langle x_{i}, x | \eta \rangle )^*\\&=\langle \sum_{i\in \mathbb{N}}\langle x , x_{i} |\xi\rangle x_{i}, x | \eta \rangle^*\\&=\langle S_{\xi} x, x | \eta \rangle^*.
\end{align*}
\end{proof}
\end{corollary}
\section{Tensor product of $\mathcal A$-2-Frames }
Let $\mathcal A$ and $\mathcal B$ be $C^*$-algebras, $E$ an $\mathcal A$-2-Hilbert space and $F$ be a $\mathcal B$-2-Hilbert space. We take $\mathcal A\otimes \mathcal B$ as the completion of $\mathcal A\otimes_{alg} \mathcal B$ with the spatial norm. Hence  $\mathcal A\otimes \mathcal B$ is a $C^*$-algebra and for every $a\in \mathcal A $ and $b\in \mathcal B$ we have $\|a\otimes b\|=\|a\|\|b\|$. The algebraic tensor product $ E\otimes_{alg} F$ is a pre-Hilbert $\mathcal A\otimes \mathcal B$-module with module action 
\begin{align*}
(a\otimes b)(x\otimes y)=ax\otimes by     \qquad (a\in A, b\in B, x\in E, y\in F )
\end{align*}
and  $\mathcal A\otimes \mathcal B$-valued 2-inner product
\begin{align*}
\langle x_{1}\otimes y_{1}, x_{2}\otimes y_{2} | \xi \otimes \eta  \rangle=\langle x_{1},  x_{2} | \xi  \rangle \otimes\langle y_{1},  y_{2} | \eta \rangle  \qquad (x_{1},  x_{2}, \xi \in E, y_{1},  y_{2}, \eta \in F).
\end{align*}
The $\mathcal A\otimes \mathcal B$-2-norm on $ E\otimes F$ is defined by
\begin{align*}
\|x_{1}\otimes x_{2}, y_{1}\otimes y_{2}\|=\|x_{1}, y_{1}\|_{\mathcal A}\|x_{2}, y_{2}\|_{\mathcal B}  \qquad (x_{1}, y_{1}\in E, x_{2}, y_{2}\in F ).
\end{align*}
where $\|. , . \|_{\mathcal A}$ and $\|. , . \|_{\mathcal B}$ are norms generated by $\langle . , . | . \rangle_{\mathcal A}$ and $\langle . , . | . \rangle_{\mathcal B}$ respectively. The space $ E\otimes F$ is complete with the above 2-inner product. Therefore, the space $ E\otimes F$ is an  $\mathcal A\otimes \mathcal B$-2-Hilbert space.\\
The following definition is the extension of (2.1) to the sequence $\{x _{i}\otimes y_{i}\}_{i\in \mathbb{N}}$.
\begin{definition}
Let $\{x _{i}\}$ and $\{y _{i}\}$ be two sequences in $\mathcal A$-2-Hilbert space $E$ and $\mathcal B$-2-Hilbert space $F$,  respectively. Then, the sequence $\{x _{i}\otimes y_{i}\}_{i\in \mathbb{N}}$ is said to be a tensor product of $\mathcal A\otimes \mathcal B$-2-frame for the tensor product of $\mathcal A\otimes \mathcal B$-2-Hilbert space $ E\otimes F$ associated to $\xi \otimes \eta$ if there exist two constants $0<A\le B<\infty$ such that
\begin{align*}
A\langle x\otimes y, x\otimes y | \xi \otimes \eta  \rangle\leq \sum_{i,j}\langle x\otimes y, x_{i}\otimes y_{j} | \xi \otimes \eta  \rangle\langle x_{i}\otimes y_{j}, x\otimes y | \xi \otimes \eta  \rangle\leq B\langle x\otimes y, x\otimes y | \xi \otimes \eta  \rangle
\end{align*} 
for all $x\in  E, y\in F$.

The numbers $A$ and $B$ are called lower and upper frame bounds of the tensor product of  $\mathcal A\otimes \mathcal B$-2-frame respectively.
\end{definition}
\begin{lemma}
Let  $\{x _{i}\}_{i\in I}$ be an $\mathcal A$-2-frame for $E$ with frame bounds $A$ and $B$ associated to $\xi$ and let $\{y _{i}\}_{j\in J}$  be a  $\mathcal B$-2-frame for $F$ with frame bounds $C$ and $D$ associated with $\eta$. Then  $\{x _{i}\otimes y_{i}\}_{i\in I, j\in J}$ is an $\mathcal A\otimes \mathcal B$-2-frame for $ E\otimes F$ with frame bounds $AC$ and $BD$ associated with $\xi\otimes \eta $,  if  $\{x _{i}\}_{i\in I}$ and $\{y _{i}\}_{j\in J}$ are tight or Parseval frames, then so is  $\{x _{i}\otimes y_{i}\}_{i\in I, j\in J}$ .
\end{lemma}
\begin{proof}
Let $x\in E$ and $y\in F$. Then we have
\begin{align} 
A\langle x, x |\xi \rangle\leq \sum_{i\in I}\langle x, x_{i} |\xi \rangle\langle x_{i}, x |\xi \rangle\leq B\langle x, x |\xi \rangle
\end{align}
\begin{align} 
C\langle y, y |\eta \rangle\leq \sum_{j\in J}\langle y, y_{i} |\eta \rangle\langle y_{i}, y |\eta \rangle\leq B\langle y, y |\eta \rangle
.\end{align} 
We know if $a$, $b$ are hermitian elements of $\mathcal A$ and $a\leq b$, then for every positive element $x$ of $\mathcal B$, we have $a\otimes x \leq b\otimes x$. Therefore
\begin{align*}
A\langle x, x |\xi \rangle\otimes \langle y, y |\eta \rangle\leq \sum_{i\in I}\langle x, x_{i} |\xi \rangle\langle x_{i}, x |\xi \rangle \otimes \langle y, y |\eta \rangle\leq B\langle x, x |\xi \rangle\otimes \langle y, y |\eta \rangle.
\end{align*}
Now by (3.2), we have
\begin{align*} 
AC\langle x, x |\xi \rangle\otimes \langle y, y |\eta \rangle&\leq \sum_{i\in I} \sum_{j\in J}\langle x, x_{i} |\xi \rangle\langle x_{i}, x |\xi \rangle \otimes \langle y, y_{i} |\eta \rangle\langle y_{i}, y |\eta \rangle \\&\leq  B\langle x, x |\xi \rangle\otimes \sum_{j\in J}\langle y, y_{i} |\eta \rangle\langle y_{i}, y |\eta \rangle\\& \leq  BD\langle x, x |\xi \rangle\otimes \langle y, y |\eta \rangle.
\end{align*}
Consequently, we have
\begin{align*} 
AC\langle x\otimes y, x\otimes y | \xi \otimes \eta  \rangle&\leq \sum_{i,j}\langle x\otimes y, x_{i}\otimes y_{j} | \xi \otimes \eta  \rangle\langle x_{i}\otimes y_{j}, x\otimes y | \xi \otimes \eta  \rangle\\&\leq BD\langle x\otimes y, x\otimes y | \xi \otimes \eta  \rangle.
\end{align*}
From these inequalities, it follows that for all $z=\sum_{k=1}^{n}u_{k}\otimes v_{k} \in E\otimes _{alg} F$ 
\begin{align*} 
AC\langle z, z | \xi \otimes \eta  \rangle&\leq \sum_{i,j}\langle z, z | \xi \otimes \eta  \rangle\langle x_{i}\otimes y_{j}, z | \xi \otimes \eta  \rangle\\&\leq BD\langle z, z | \xi \otimes \eta  \rangle.
\end{align*} 
Hence it  holds for all $z\in E\otimes  F$.
\end{proof}
If take $\mathcal A=\mathbb{C}$ and define $\mathbb{C}\otimes \mathcal B$-2-product by 
\begin{align}\label{3.3}
\langle x_{1}\otimes y_{1}, x_{2}\otimes y_{2} | \xi \otimes \eta  \rangle=\langle x_{1},  x_{2} | \xi  \rangle \langle y_{1},  y_{2} | \eta \rangle  \qquad (x_{1},  x_{2}, \xi \in E, y_{1},  y_{2}, \eta \in F).
\end{align}
We have the following proposition.
\begin{proposition}
Let $\{x _{i}\otimes y_{i}\}_{i\in I, j\in J}$ be a $\mathbb{C}\otimes \mathcal B$-2-frame for $ E\otimes F$ associated with $\xi\otimes \eta $, then $\{y _{i}\}_{j\in J}$  is a  $\mathcal B$-2-frame for $F$ associated with $\eta$. 
\begin{proof}
Suppose that $\{x _{i}\otimes y_{i}\}_{i\in I, j\in J}$ is a $\mathbb{C}\otimes \mathcal B$-2-frame for $ E\otimes F$ associated to $\xi\otimes \eta $. Then for $x\otimes y\in  E\otimes F-\{0\otimes0\}$, we have
\begin{align*}
A\langle x\otimes y, x\otimes y | \xi \otimes \eta  \rangle\leq \sum_{i,j}\langle x\otimes y, x_{i}\otimes y_{j} | \xi \otimes \eta  \rangle\langle x_{i}\otimes y_{j}, x\otimes y | \xi \otimes \eta  \rangle\leq B\langle x\otimes y, x\otimes y | \xi \otimes \eta  \rangle.
\end{align*} 
By using (\ref{3.3}) the above equation becomes
\begin{align*} 
A\langle x, x |\xi \rangle \langle y, y |\eta \rangle&\leq \sum_{i\in I} |\langle x, x_{i} |\xi \rangle|^2  \sum_{j\in J}\langle y, y_{j} |\eta \rangle\langle y_{j}, y |\eta \rangle \\& \leq  B\langle x, x |\xi \rangle \langle y, y |\eta \rangle.
\end{align*}
This gives
\begin{align*} 
\frac {A\langle x, x |\xi \rangle}{ \sum_{i\in I} |\langle x, x_{i} |\xi \rangle|^2} \langle y, y |\eta \rangle&\leq  \sum_{j\in J}\langle y, y_{j} |\eta \rangle\langle y_{j}, y |\eta \rangle \\& \leq \frac{ B\langle x, x |\xi \rangle}{ \sum_{i\in I} |\langle x, x_{i} |\xi \rangle|^2} \langle y, y |\eta \rangle.
\end{align*}
Therefore 
\begin{align*} 
A_{1} \langle y, y |\eta \rangle \leq  \sum_{j\in J}\langle y, y_{j} |\eta \rangle\langle y_{j}, y |\eta \rangle \leq B_{1} \langle y, y |\eta \rangle. \qquad (\forall y\in F)
\end{align*}
Where $A_{1}=\frac {A\langle x, x |\xi \rangle}{ \sum_{i\in I} |\langle x, x_{i} |\xi \rangle|^2}$ and $B_{1}= \frac{ B\langle x, x |\xi \rangle}{ \sum_{i\in I} |\langle x, x_{i} |\xi \rangle|^2}.$
\end{proof}
\end{proposition}
\begin{remark} 
If the sequences $\{x_{i}\}_{i\in I}$,  $\{y_{j}\}_{j\in J} $ and$ \{x_{i}\otimes y_{j}\}_{i\in I,j\in J}$ are 2-frames for Hilbert spaces $E_{\xi}$, $F_{\eta}$ and $(E\otimes F)_{\xi \otimes \eta}$ respectively and $S_{\xi}$, $S_{\eta}$ and $S_{xi \otimes \eta}$ are frame operators of above frames respectively, then from (\ref{2.2}), we have the following
\begin{align*} 
S_{\xi}x=\sum_{i\in I}\langle x , x_{i} |\xi\rangle x_{i},
S_{\eta}y=\sum_{j\in J}\langle y , y_{j} |\eta\rangle y_{j}.
\end{align*}
\begin{align*}
S_{\xi\otimes \eta}(x\otimes y)=\sum_{i\in I ,j\in J}\langle x\otimes y , x_{i}\otimes y_{j} |\xi \otimes \eta\rangle( x_{i}\otimes y_{j}),( x\in E, y\in F, x\otimes y\in E\otimes F )
\end{align*}
\end{remark}
\begin{theorem}
If the sequences $\{x_{i}\}_{i\in I}$,  $\{y_{j}\}_{j\in J} $ and$ \{x_{i}\otimes y_{j}\}_{i\in I,j\in J}$ are 2-frames for Hilbert spaces $E_{\xi}$, $F_{\eta}$ and $(E\otimes F)_{\xi \otimes \eta}$, respectively and $S_{\xi}$, $S_{\eta}$ and $S_{\xi \otimes \eta}$ are frame operators respectively, then $S_{\xi \otimes \eta}=S_{\xi}\otimes S_{\eta}$.
\begin{proof}
For $x\otimes y\in E\otimes F$, we have
\begin{align*}
S_{\xi\otimes \eta}(x\otimes y)&=\sum_{i\in I ,j\in J}\langle x\otimes y , x_{i}\otimes y_{j} |\xi \otimes \eta\rangle( x_{i}\otimes y_{j})\\&=\sum_{i\in I ,j\in J}\langle x, x_{i} | \xi \rangle_{\mathcal A} \langle y , y_{j} |\eta \rangle_{\mathcal B}(x_{i}\otimes y_{j})\\&=\sum_{i\in I}\langle x, x_{i} | \xi \rangle_{\mathcal A}x_{i}\otimes \sum_{j\in J}\langle y, y_{j} | \eta \rangle_{\mathcal B}y_{i}\\&=S_{\xi}x\otimes S_{\eta}y=(S_{\xi}\otimes S_{\eta})(x\otimes y).
\end{align*}
Hence $S_{\xi \otimes \eta}=S_{\xi}\otimes S_{\eta}$.
\end{proof}
\end{theorem}
\begin{lemma}
Suppose that $\{x_{i}\}_{i\in I}$ is a sequence in $\mathcal A$-2-Hilbert space $E$, with $x=\sum_{i\in I}\langle x, x_{i} | \xi \rangle x_{i}$ holds for all $x\in E$, then  $\{x_{i}\}_{i\in I}$ is an $\mathcal A$-2-normalized tight frame for $E$.
\begin{proof}
$x=\sum_{i\in I}\langle x, x_{i} | \xi \rangle x_{i}$, so $\langle x , x | \xi \rangle=\langle \sum_{i\in I}\langle x, x_{i} | \xi \rangle x_{i} , x | \xi \rangle$, hence $\Longrightarrow \langle x , x | \xi \rangle=\sum_{i\in I}\langle x, x_{i} | \xi \rangle\langle x_{i} , x | \xi \rangle.$
\end{proof}
\end{lemma}
\begin{corollary}
Assume that  $\{x_{i}\otimes y_{j}\}_{i\in I, j\in J}$ is a sequence in $\mathcal A\otimes \mathcal B$-2-Hilbert space $E\otimes F$ and $x\otimes y=\sum_{i\in I ,j\in J}\langle x\otimes y , x_{i}\otimes y_{j} |\xi \otimes \eta\rangle( x_{i}\otimes y_{j})$, for all $x\in E, y\in F$. Then $\{x_{i}\otimes y_{j}\}_{i\in I, j\in J}$ is an $\mathcal A\otimes \mathcal B$-2-normalized tight frame for $E\otimes F$.
\begin{proof}
By lemma (3.6) is clear.
\end{proof}
\end{corollary}

\end{document}